\documentclass[a4paper,11pt]{article}
\usepackage[latin1]{inputenc}
\usepackage[T1]{fontenc}
\usepackage[francais]{babel}
\usepackage{amssymb}
\usepackage{amsfonts, amsmath}
 \newtheorem{thm}{THÉORÈME}
 \newtheorem{cor}{Corollaire}

 \newtheorem{rem}{REMARQUE}
 \newtheorem{Exmp}{Exemple}

\title{\textbf{Le th\'{e}or\`{e}me de Riemann-Roch et ses applications}}
\author{\textbf{A. Lesfari}
\\Département de Mathématiques
\\Faculté des Sciences
\\Université Chouaïb Doukkali
\\B.P. 20, El-Jadida, Maroc.
\\E-mail : Lesfariahmed@yahoo.fr , lesfari@ucd.ac.ma }

\date{}

\begin{document}
\maketitle Abstract. The Riemann-Roch theorem is of utmost
importance in the algebraic geometric theory of compact Riemann
surfaces. It tells us how many linearly independent meromorphic
functions there are having certain restrictions on their poles.
The aim of this article is to present a simple direct proof of
this theorem and explore some of its numerous consequences. We
also give an analytic proof of the Riemann-Hurwitz formula. As an
application, we compute the genus of some interesting algebraic
curves.

\vskip0.4cm \null\hskip0.4cm

Dans ce travail, nous allons étudier un des théorèmes les plus
importants de la théorie des surfaces de Riemann compactes : le
théorème de Riemann-Roch. Il s'agit d'un théorème d'existence
efficace qui permet, entre autres, de déterminer le nombre de
fonctions méromorphes linéairement indépendantes ayant certaines
restrictions sur leurs pôles. Le but de cet article est de donner
une preuve élémentaire, constructive bien qu'un peu technique du
théorème de Riemann-Roch. Nous mentionnons quelques conséquences
de ce théorème et nous donnons également une preuve analytique de
la formule de Riemann-Hurwitz. Elle exprime le genre d'une surface
de Riemann à l'aide du nombre de ses points de ramifications et du
nombre de ses feuillets. Nous montrons que cette formule fournit
un moyen efficace pour déterminer le genre d'une surface de
Riemann donnée. Plusieurs exemples intéressants seront étudiés.

\section{Préliminaires}

Nous commençons par des généralités sur les surfaces de Riemann
compactes $X$ ou courbes algébriques complexes. Ce sont des
variétés analytiques de dimension $1$ complexe ($2$ réelle) munies
d'atlas dont les changements de cartes sont holomorphes. On les
définit par une équation de la forme
$$F(w,z) \equiv p_{0}(z)w^{n}+p_{1}(z)w^{n-1}+\cdots +p_{n}(z)=0,$$
où $F(w,z)$ est un polynôme à deux variables complexes $w$ et $z$,
de degré $n$ en $w$ et irréductible. Ici $p_{0}(z)\neq 0,$
$p_{1}(z),\ldots ,$ $p_{n}(z)$ sont des polynômes en $z.$ La
surface $X$ est homéomorphe à un tore à $g$ trous. Le nombre $g$
s'appelle genre de la surface $X$.

Soient $p$ un point de $X$, $\tau_p:X\rightarrow
\overline{\mathbb{C}}$ un paramètre local en $p$ (ou une
uniformisante locale en $p$, i.e., une carte locale en $p$
appliquant $p$ sur $0$) et $f$ une fonction méromorphe au
voisinage de $p$. L'ordre de $f$ en $p$, est l'unique entier $n$
tel que : $f=\tau_p^n.g,$ où $g$ est holomorphe ne s'annulant pas
en $p$. Dans le cas où $f=0$, on choisit par convention
$n=+\infty$. L'entier $n$ dépend de $p$ et de $f$ et on le note
$\mbox{ord}_pf$. On a $$\mbox{ord}_p(f+g)\geq \inf (\mbox{ord}_pf
, \mbox{ord}_p g),$$ et $$\mbox{ord}_p
fg=\mbox{ord}_pf+\mbox{ord}_pg.$$

Un diviseur sur une surface de Riemann $X$ est une combinaison
formelle du type
$$\mathcal{D}=\sum_{p\in X}n_p.p=\sum_j n_{j}p_{j},\quad n_j\in
\mathbb{Z},$$ avec $(p_{j})$ une famille localement finie de
points de $X$. Le diviseur $\mathcal{D}$ est fini si son support
est fini et ce sera toujours le cas si $X$ est une surface de
Riemann compacte. L'ensemble des diviseurs sur $X$ est un groupe
abélien noté $\mbox{Div} (X)$. L'addition des diviseurs est
définie par l'addition des coefficients. Soit $f\neq 0$, une
fonction méromorphe sur $X$. A cette fonction $f$, on fait
correspondre un diviseur noté $(f)$ en posant $$(f)=\sum_{p\in
X}\mbox{ord}_pf . p ,$$ où les $\mbox{ord}_pf$ sont nuls sauf un
nombre fini d'entre eux. En désignant par $\alpha_1,...,\alpha_l$
les zéros de $f$ de multiplicité $n_1,...,n_l$ respectivement et
par $\beta_1,...,\beta_m$ les pôles de $f$ de multiplicité
$p_1,...,p_m$ respectivement, on obtient
\begin{eqnarray}
(f)&=&\sum_{j=1}^l n_{j}\alpha_{j}-\sum_{j=1}^m p_{j}\beta_{j},\nonumber\\
&=&(\text{diviseur des zéros de }f)-(\text{diviseur des pôles
de}f).\nonumber
\end{eqnarray}
On a
\begin{eqnarray}
(fg)&=&(f)+(g),\nonumber\\
(f^{-1})&=&-(f),\nonumber\\
(f)&=&(g)\Longrightarrow \frac{f}{g}=\mbox{constante}.\nonumber
\end{eqnarray}
Tout diviseur d'une fonction méromorphe est dit diviseur
principal. L'ensemble des diviseurs principaux forme un
sous-groupe $\mbox{Div}^\circ (X)$ de $\mbox{Div} (X)$. Le groupe
quotient $\mbox{Div} (X)/\mbox{Div}^\circ (X)=\mbox{Pic} (X),$ est
le groupe de Picard de $X$. Le degré du diviseur $\mathcal{D}$ est
un entier noté $\mbox{deg} \mathcal{D}$ et est défini par
$\mbox{deg} \mathcal{D}=\sum_jn_j$. On dit qu'un diviseur
$\mathcal{D}$ est positif (ou effectif) et on note
$\mathcal{D}\geq 0$, si les entiers $n_{j}$ qui interviennent dans
la somme sont positifs. Un diviseur $\mathcal{D}_1$ est plus grand
qu'un diviseur $\mathcal{D}_2$ si $\mathcal{D}_1-\mathcal{D}_2$
est positif. Deux diviseurs $\mathcal{D}_1$ et $\mathcal{D}_2$
sont dits linéairement équivalents (et on note $\mathcal{D}_1\sim
\mathcal{D}_2$) si $\mathcal{D}_1-\mathcal{D}_2$ est principal,
i.e., si $\mathcal{D}_1-\mathcal{D}_2=(f)$ où $f$ est une fonction
méromorphe. L'application $$\mbox{deg} : \mbox{Div}
(X)\longrightarrow \mathbb{Z},\quad
\mathcal{D}\longmapsto\mbox{deg} (\mathcal{D}),$$ est un
homorphisme de groupe. Sur toute surface de Riemann compacte, une
fonction méromorphe $f\neq 0$ a le même nombre de zéros que des
pôles, donc $\mbox{deg}(f)=0$. Autrement dit, tout diviseur
principal a le degré $0$.

Si $\mathcal{D}=\sum_{p\in X}n_p . p$ est un diviseur, on notera
$\mathcal{L}(\mathcal{D})$ l'ensemble des fonctions méromorphes
$f$ telles que: $\mbox{ord}_pf+n_p\geq 0$, pour tout $p\in X$.
Autrement dit,
$$\mathcal{L}(\mathcal{D})=\{f\text{méromorphe sur }X:(f)+\mathcal{D}\geq 0\},$$
i.e., l'espace vectoriel des fonctions méromorphes $f$ telles que:
$(f)+\mathcal{D}\geq 0.$ Si $(f)+\mathcal{D}$ n'est $\geq 0$ pour
aucun $f,$ on posera $\mathcal{L}(\mathcal{D})=0$. Par exemple, si
le diviseur $\mathcal{D}$ est positif alors
$\mathcal{L}(\mathcal{D})$ est l'ensemble des fonctions
holomorphes en dehors de $\mathcal{D}$ et ayant au plus des
p\^{o}les simples le long de $\mathcal{D}$. Notons que si
$\mathcal{D}_1\sim \mathcal{D}_2$, alors
$\mathcal{L}(\mathcal{D}_1)$ est isomorphe à
$\mathcal{L}(\mathcal{D}_2)$, d'où $\dim
\mathcal{L}(\mathcal{D}_1)=\dim \mathcal{L}(\mathcal{D}_2)$.

Une forme différentielle sur $X$ s'\'{e}crit $\omega =f(\tau)
d\tau$, où $\tau$ est le paramètre local et $f$ une fonction
complexe de $\tau .$ On dit que $\omega $ est une différentielle
abélienne si $f(\tau)$ est une fonction méromorphe sur $X$,
holomorphe si $f(\tau)$ est une fonction holomorphe sur $X$, ayant
un pôle d'ordre $k$ ou un zéro d'ordre $k$ en un point $p$ si
$f(\tau)$ a un pôle d'ordre $k$ ou un zéro d'ordre $k$ en ce
point. L'ensemble des différentielles holomorphes sur
$\mathcal{C}$ est de dimension $g.$ On peut associer à chaque
forme différentielle $\omega$ un diviseur noté $(\omega)$. Si
$\omega=fd\tau_p$ avec $f$ une fonction méromorphe sur $X$ et
$\tau_p$ un paramètre local en $p\in X$, on définit l'ordre de
$\omega$ en $p$ par $\mbox{ord}_p\omega=\mbox{ord}_0f$ et le
diviseur $(\omega)$ de $\omega$ par $(\omega)=\sum_{p\in
X}\mbox{ord}_p\omega.p$.

Si $\mathcal{D}=\sum_{p\in X}n_p . p$ est un diviseur, on définit
de façon analogue à $\mathcal{L}(\mathcal{D})$, l'espace linéaire
$\Omega(\mathcal{D})$ comme étant l'ensemble des formes
différentielles méromorphes $\omega$ sur $X$ telles que :
$\mbox{ord}_p\omega+n_p\geq 0$, pour tout $p\in X$. Autrement dit,
$$\Omega(\mathcal{D})=\{\omega\text{méromorphe sur }X:(\omega)+\mathcal{D}\geq 0\},$$
c'est-à-dire l'ensemble des formes différentielles méromorphes
$\omega$ sur $X$ telles que : $(\omega)+\mathcal{D}\geq 0$. Si
$\mathcal{D}_1\sim \mathcal{D}_2$, alors $\Omega(\mathcal{D}_1)$
est isomorphe à $\Omega(\mathcal{D}_2)$, d'où $\dim
\Omega(\mathcal{D}_1)=\dim \Omega(\mathcal{D}_2)$. Dans le cas où
le diviseur $\mathcal{D}$ est négatif, alors
$(\omega)+\mathcal{D}\geq 0$ est l'ensemble des formes
différentielles qui n'ont pas de pôles et qui ont des zéros au
moins aux points de $\mathcal{D}$. Notons enfin qu'en vertu du
théorème des résidus, on a $\sum_{p\in X}\mbox{Res} (\omega)=0$,
où $\mbox{Res} (\omega)$ est le résidu en $p$ de $\omega$, i.e.,
le coefficient de $\frac{1}{\tau_p}$ dans le développement de $f$
en série de Laurent.

Soit $\mathcal{D}$ un diviseur sur une surface de Riemann compacte
$X$ et $K$ un diviseur canonique sur $X$, i.e., le diviseur
$(\omega)$ d'une $1-$forme méromorphe $\omega\neq 0$ sur $X$.
L'application
\begin{equation}\label{eqn:euler}
\psi:\mathcal{L}(K-\mathcal{D})\longrightarrow
\Omega(-\mathcal{D}),\quad f\longmapsto f\omega,
\end{equation}
est un isomorphisme. En effet, on a
\begin{eqnarray}
(f\omega)=(f)+(\omega)&=&(f)+K,\nonumber\\
&\geq&-(K-\mathcal{D})+K,\nonumber\\
&=&-(-\mathcal{D}),\nonumber
\end{eqnarray}
ce qui montre que l'application $\psi$ est bien définie. Cette
dernière est injective, i.e., l'équation $f\omega=g\omega$
entraine $f=g$. Montrons que $\psi$ est surjective. Soit $\eta \in
\Omega(-\mathcal{D})$, d'où il existe une fonction méromorphe $h$
sur $X$ telle que : $h\omega=\eta$. Dès lors,
\begin{eqnarray}
(h)+K&=&(h)+(\omega),\nonumber\\
&=&(h\omega),\nonumber\\
&=&(\eta)\geq-(-\mathcal{D}),\nonumber
\end{eqnarray}
d'où $(h)\geq -(K-\mathcal{D})$, $h\in \mathcal{L}(K-\mathcal{D})$
et par conséquent $\psi$ est surjective.

\section{Théorème de Riemann-Roch et formule de Riemann-Hurwitz}

\begin{thm}
(de Riemann-Roch) : Si $X$ une surface de Riemann compacte et
$\mathcal{D}$ est un diviseur sur $X$, alors
\begin{equation}\label{eqn:euler}
\dim \mathcal{L}(\mathcal{D})-\dim
\mathcal{L}(K-\mathcal{D})=\mbox{deg} \mathcal{D}-g+1,
\end{equation}
où $K$ est le diviseur canonique sur $X$ et $g$ est le genre de
$X$. La formule (2) peut s'écrire sous la forme équivalente
\begin{equation}\label{eqn:euler}
\dim \mathcal{L}(\mathcal{D})-\dim \Omega(-\mathcal{D})=\mbox{deg}
\mathcal{D}-g+1.
\end{equation}
\end{thm}
\textbf{Démonstration}: L'équivalence entre les formules $(2)$ et
$(3)$ résulte immédiatement de l'isomorphisme $(1)$. La preuve du
théorème est immédiate dans le cas où $\mathcal{D}=0$ car
$\mathcal{L}(0)$ est l'ensemble des fonctions holomorphes sur $X$.
Or toute fonction holomorphe sur une surface de Riemann compacte
est constante, donc $\mathcal{L}(0)=\mathbb{C}$. En outre, on sait
que la dimension de l'espace des formes différentielles
holomorphes sur $X$ est le genre $g$ de $X$, d'où le résultat. La
preuve du théorème va se faire en plusieurs étapes :\\
\null\hskip0.4cm \emph{Étape 1}: Soit $\mathcal{D}$ un diviseur
positif. Autrement dit,
$$\mathcal{D}=\sum_{k=1}^mm_kp_k,\quad n_k\in \mathbb{N}^*,\quad
p_k\in X.$$  Soit $f\in \mathcal{L}(\mathcal{D}).$ Au voisinage de
$p_k$, on a
$$df=(\sum_{j=-n_k-1}^\infty c_j^k\tau^j)d\tau,$$
donc $df$ est méromorphe. Plus précisement, $df$ a un pôle d'ordre
$n_k+1$ en $p_k$. Comme $f$ est méromorphe, alors $df$ ne peut pas
avoir de pôle simple et dès lors son résidu est nul, i.e.,
$c_{-1}^k=0.$ Soit $\left( a_{1},\ldots ,a_{g},b_{1},\ldots
,b_{g}\right) $ une base de cycles dans le groupe d'homologie
$H_{1}\left( X,\mathbb{Z}\right) $ de telle fa\c{c}on que les
indices d'intersection de cycles deux \`{a} deux s'\'{e}crivent :
$$\left( a_{i},a_{i}\right) =\left( b_{i},b_{i}\right)
=0,\quad\left( a_{i},b_{j}\right) =\delta _{ij},\quad1\leq i,j\leq
g.$$  Posons
\begin{equation}\label{eqn:euler}
\eta=df-\sum_{k=1}^m\sum_{j=2}^{n_k}c_j^k\eta_k^j.
\end{equation}
D'où
$$\int_{a_i}\eta=\int_{a_i}df-\sum_{k=1}^m\sum_{j=2}^{n_k}c_j^k\int_{a_i}\eta_k^j.$$
L'intégrale d'une différentielle exacte le long d'un chemin fermé
étant nulle, donc $\int_{a_i}df=0$. Rappelons (analyse harmonique)
que si $X$ est une surface de Riemann de genre $g$, alors pour
tout $n\geq 2$ et pour tout $p\in X$, il existe une différentielle
holomorphe $\eta$ sur $X\backslash p$ telle que :
$$\eta=(\frac{1}{\tau^n}+\circ (\tau))d\tau,$$ et en outre, on a
$\int_{a_i}\eta=0$. La forme $\eta$ étant holomorphe, alors
$$\eta=c_1\omega_1+...+c_g\omega_g,$$ où $(\omega_1,...,\omega_g)$
est une base de $\Omega(X)$ et dès lors
$$\eta=c_1\int_{a_i}\omega_1+...+c_g\int_{a_i}\omega_g,\quad i=1,...,g$$
Puisque la matrice $$E=(\int_{a_i}\omega_j)_{1\leq i,j\leq g},$$
est inversible, alors $c_1=...=c_g=0$, donc $\eta=0$ et d'après
(4), on a
$$df=\sum_{k=1}^m\sum_{j=2}^{n_k}c_j^k\eta_k^j.$$
Considèrons l'application
$$\varphi:\mathcal{L}(\mathcal{D})\longrightarrow
V\equiv\{(c_j^k):\sum_{k=1}^m\sum_{j=2}^{n_k}c_j^k\int_{b_l}\eta_k^j=0\},\quad
f\longmapsto c_j^k.$$ Notons que
\begin{eqnarray}
\mbox{Ker} \varphi&=&\{f : \mbox{méromorphe sur X et n'ayant pas
de
pôle}\},\nonumber\\
&=&\{f : f\mbox{est une constante}\},\nonumber
\end{eqnarray}
d'où $\dim \mbox{Ker} \varphi=1$ et par conséquent, $$\dim
\mathcal{L}(\mathcal{D})=\dim V+1.$$ Les espaces
$\frac{\mathcal{L}(\mathcal{D})}{\mathbb{C}}$ et $V$ sont
isomorphes et on a
\begin{eqnarray}
\dim\mathcal{L}(\mathcal{D})-1&=&\dim V,\nonumber\\
&=&\dim\{(c_j^k):\sum_{k=1}^m\sum_{j=2}^{n_k}c_j^k\int_{b_l}\eta_k^j=0\},\nonumber\\
&=&\mbox{deg}\mathcal{D}-\mbox{rang}\mathcal{M},\nonumber
\end{eqnarray}
où
$$
{\mathcal{M}}=\left(\begin{array}{ccccccccc}
\int_{b_l}\eta_1^2&\int_{b_l}\eta_1^3&...&\int_{b_l}\eta_1^{n_1+1}&\int_{b_l}\eta_2^2&...&
\int_{b_l}\eta_2^{n_2+1}&...&\int_{b_l}\eta_m^{n_m+1}\\
\int_{b_2}\eta_1^2&\int_{b_2}\eta_1^3&...&\int_{b_2}\eta_1^{n_1+1}&\int_{b_2}\eta_2^2&...&
\int_{b_2}\eta_2^{n_2+1}&...&\int_{b_2}\eta_m^{n_m+1}\\
\vdots&\vdots&...&\vdots&\vdots&...&\vdots&...&\vdots\\
\int_{b_g}\eta_1^2&\int_{b_g}\eta_1^3&...&\int_{b_g}\eta_1^{n_1+1}&\int_{b_g}\eta_2^2&...&
\int_{b_g}\eta_2^{n_2+1}&...&\int_{b_g}\eta_m^{n_m+1}
\end{array}\right),
$$
est la matrice dont le nombre de lignes est $g$ et le nombre de
colonnes est $\mbox{deg} \mathcal{D}$. Notons que
\begin{eqnarray}
\mbox{rang} \mathcal{M}&=&\mbox{Nombre de colonnes}-
\mbox{Nombre de relations entre ces colonnes},\nonumber\\
&=&\mbox{deg} \mathcal{D}-\dim V,\nonumber\\
&=&\mbox{deg} \mathcal{D}-\dim \mathcal{L}(\mathcal{D})+1.
\end{eqnarray}
Calculons maintenant le rang de $\mathcal{M}$ d'une autre façon.
Soit $(\omega_1,...,\omega_g)$ une base normalisée de $\Omega(X)$,
i.e., de telle sorte que : $\int_{a_i}\omega_j=\delta_{ij}$. Au
voisinage de $p_k$, la forme $\omega_s$ admet un dévelopement en
série de Taylor,
$$w_s=(\sum_{j=0}^\infty \alpha_{sj}^k\tau^j)d\tau.$$
Posons $\varphi_s\equiv \int_0^z\omega_s$ et soit $X^*$ la
représentation normale de la surface de Riemann $X$, i.e., un
polygône à $4g$ côtés identifiés deux à deux (On le désigne par
$(a_1b_1a_1^{-1}b_1^{-1}...a_gb_ga_g^{-1}b_g^{-1})$ et peut être
définit à partir d'une triangulation de la surface $X$). Notons
que si $\tau\in a_j$, alors il est identifié à $\tau^*\in
a_j^{-1}$, d'où
$$\varphi_s(\tau^*)=\varphi_s(\tau)+\int_{b_j}\omega_s.$$ De même,
si $\tau\in b_j$, alors il est identifié à $\tau^*\in b_j^{-1}$ et
$$\varphi_s(\tau^*)=\varphi_s(\tau)+\int_{a_j}\omega_s.$$ On a
\begin{eqnarray}
\int_{\partial X^*}\varphi_s\eta_k^n
&=&\sum_{j=1}^g(\int_{a_j}\varphi_s\eta_k^n+\int_{b_j}\varphi_s\eta_k^n
+\int_{a_j^{-1}}(\varphi_s+\int_{b_j}\omega_s)\eta_k^n
+\int_{b_j^{-1}}(\varphi_s-\int_{a_j}\omega_s)\eta_k^n),\nonumber\\
&=&\sum_{j=1}^g(-\int_{b_j}\omega_s\int_{a_j}\eta_k^n
+\int_{a_j}\omega_s\int_{b_j}\eta_k^n),\nonumber\\
&=&\sum_{j=1}^g(-\omega_s(b_j)\eta_k^n(a_j)+\omega_s(a_j)\eta_k^n(b_j)),\nonumber\\
&=&\sum_{j=1}^g\omega_s(a_j)\eta_k^n(b_j),\nonumber\\
&=&\eta_k^n(b_s).
\end{eqnarray}
Or
\begin{eqnarray}
\int_{\partial X^*}\varphi_s\eta_k^n&=& 2\pi i\sum_k
\mbox{Rés}_{p_k}(\varphi_s\eta_k^n),\nonumber\\
&=&2\pi i \frac{\alpha_{s,n-2}^k}{n-1},\nonumber
\end{eqnarray}
donc la matrice $\mathcal{M}$ a comme coefficient
$$\int_{b_s}\eta_k^n=\eta_k^n(b_s)=2\pi i \frac{\alpha_{s,n-2}^k}{n-1}.$$
Dès lors $$\det \mathcal{M}=C\det \mathcal{N},$$ où $$C\equiv(2\pi
i)(\pi i)...(\frac{2\pi i}{n_1})(2\pi i)...(\frac{2\pi
i}{n_2})...(\frac{2\pi i}{n_m}),$$ est une constante et
$$
{\mathcal{N}}= \left(\begin{array}{ccccccccc}
\alpha_{1,0}^1&\alpha_{1,1}^1&...&\alpha_{1,n_1-1}^1&\alpha_{1,0}^2&...
&\alpha_{1,n_2-1}^2&...&\alpha_{1,n_m-2}^m\\
\alpha_{2,0}^1&\alpha_{2,1}^1&...&\alpha_{2,n_1-1}^1&\alpha_{2,0}^2&...
&\alpha_{2,n_2-1}^2&...&\alpha_{2,n_m-2}^m\\
\vdots&\vdots&...&\vdots&\vdots&...&\vdots&...&\vdots\\
\alpha_{g,0}^1&\alpha_{g,1}^1&...&\alpha_{g,n_1-1}^1&\alpha_{g,0}^2&...
&\alpha_{g,n_2-1}^2&...&\alpha_{g,n_m-2}^m\\
\end{array}\right).
$$
Calculons maintenant la dimension de l'espace
$\mathcal{L}(K-\mathcal{D})$ ou ce qui revient au même de l'espace
$\Omega(-\mathcal{D})$, i.e., celui des formes différentielles
$\omega$ qui s'annulent $n_k$ fois au point $p_k$. On a
$$\omega=\sum_{s=1}^gX_s\omega_s=
\sum_{s=1}^gX_s(\alpha_{s,0}^k+\alpha_{s,1}^k\tau+\alpha_{s,2}^k\tau^2+...)d\tau.$$
Pour que $\omega$ s'annule $n_k$ fois au point $p_k$, il faut que
les $n_k$ premières termes dans l'expression ci-dessus soient
nulles. Dès lors,
$$(X_1,...,X_g).\mathcal{N}=0,$$
tandis que la dimension de $\mathcal{L}(K-\mathcal{D})$ coincide
avec celle de l'ensemble de $(X_1,...,X_g)$ tel que :
$\omega=\sum_{s=1}^gX_s\omega_s$ s'annule $n_k$ fois au point
$p_k$, i.e.,
\begin{eqnarray}
\dim \mathcal{L}(K-\mathcal{D})&=&g-\mbox{rang}
\mathcal{N},\nonumber\\
&=&g-\mbox{rang} \mathcal{M}.\nonumber
\end{eqnarray}
D'où
$$\mbox{rang} \mathcal{M}=g-\dim \mathcal{L}(K-\mathcal{D}),$$
et en tenant compte de (5), on obtient finalement
$$\dim \mathcal{L}(\mathcal{D})-\dim \mathcal{L}(K-\mathcal{D})=\mbox{deg}
\mathcal{D}-g+1.$$ \null\hskip0.4cm \emph{Étape 2}: La preuve
donnée dans l'étape 1 est valable pour tout diviseur linéairement
équivalent à un diviseur positif étant donné que $\dim
\mathcal{L}(\mathcal{D}),$ $\dim \mathcal{L}(K-\mathcal{D})$ (ou
$\dim \Omega(-\mathcal{D})$)
et $\mbox{deg} \mathcal{D}$ ne seront pas affectés.\\
\null\hskip0.4cm \emph{Étape 3}: Soit $f$ une fonction méromorphe,
$\mathcal{D}$ un diviseur positif et posons
$\mathcal{D}'=(f)+\mathcal{D}_0,$ autrement dit, $\mathcal{D}'$ et
$\mathcal{D}_0$ sont linéairement équivalent. Nous allons tout
d'abord démontrer les assertions suivantes :
\begin{eqnarray}
&\textbf{(i)}&\dim \mathcal{L}(\mathcal{D}')=\dim
\mathcal{L}(\mathcal{D}_0).\nonumber\\
&\textbf{(ii)}&\dim \mathcal{L}(K-\mathcal{D'})=\dim
\mathcal{L}(K-\mathcal{D}_0).\nonumber\\
&\textbf{(iii)}&\mbox{deg} \mathcal{D}'=\mbox{deg}
\mathcal{D}_0.\nonumber
\end{eqnarray}
En effet, pour tout $g\in \dim \mathcal{L}(\mathcal{D}')$, on a
$$(g)+\dim \mathcal{L}(\mathcal{D}')\geq 0,$$ et dès lors
\begin{eqnarray}
(fg)+\dim
\mathcal{L}(\mathcal{D}_0)&=&(f)+(g)+\dim\mathcal{L}(\mathcal{D}_0),\nonumber\\
&=&\dim\mathcal{L}(\mathcal{D}')-\dim\mathcal{L}(\mathcal{D}_0)+(g)+\dim
\mathcal{L}(\mathcal{D}_0),\nonumber\\
&=&\dim \mathcal{L}(\mathcal{D}')+(g)\geq 0.\nonumber
\end{eqnarray}
L'application $$\mathcal{L}(\mathcal{D}')\longrightarrow
\mathcal{L}(\mathcal{D}_0),\quad g\longmapsto fg,$$ est linéaire
et admet comme réciproque
$$\mathcal{L}(\mathcal{D}_0)\longrightarrow
\mathcal{L}(\mathcal{D}'), \quad g\longmapsto \frac{g}{f}.$$ D'où
$\mathcal{L}(\mathcal{D}')\cong \mathcal{L}(\mathcal{D}_0)$ et par
conséquent
$\dim\mathcal{L}(\mathcal{D}')=\dim\mathcal{L}(\mathcal{D}_0)$. En
ce qui concerne $(ii)$, il suffit d'utiliser un raisonnement
similaire au précédent. Pour $(iii)$, on a
$$\mbox{deg}(f)=\mbox{deg}\mathcal{D}'-\mbox{deg} \mathcal{D}_0,$$ et
le résultat découle du fait que tout diviseur principal a le degré
$0$. En visageons maintenant les différents cas possibles :\\
$1^{\mbox{re}} cas$ : $\dim\mathcal{L}(\mathcal{D})>0$. Soit
$f_0\in \mathcal{L}(\mathcal{D})$, d'où
$$\mathcal{L}(\mathcal{D})\equiv (f_0)+\mathcal{D}>0,$$ et
$$\dim\mathcal{L}((f_0)+\mathcal{D})-\dim\mathcal{L}(K-(f_0)-\mathcal{D})
=\mbox{deg} ((f_0)+\mathcal{D})-g+1,$$ i.e.,
$$\dim\mathcal{L}(\mathcal{D})-\dim\mathcal{L}(K-\mathcal{D})
=\mbox{deg} \mathcal{D}-g+1.$$ $2^{\mbox{ème}} cas$ :
$\dim\mathcal{L}(\mathcal{D})=0$ et
$\dim\mathcal{L}(K-\mathcal{D})\neq 0$. En appliquant la formule
ci-dessus à $K-\mathcal{D}$, on obtient
\begin{equation}\label{eqn:euler}
\dim\mathcal{L}(K-\mathcal{D})-\dim\mathcal{L}(\mathcal{D})
=\mbox{deg} (K-\mathcal{D})-g+1.
\end{equation}
Pour la suite, on aura besoin du résultat intéressant suivant :
Pour tout diviseur canonique $K$ sur une surface de Riemann
compacte $X$, on a
\begin{equation}\label{eqn:euler}
\mbox{deg} K=2g-2.
\end{equation}
où $g$ est le genre de $X$. En effet, en posant $\mathcal{D}=K$
dans la formule (2), on obtient $$\dim
\mathcal{L}(\mathcal{D})-\dim
\mathcal{L}(K-\mathcal{D})=\mbox{deg} \mathcal{D}-g+1.$$ Or
$\mathcal{L}(0)=\mathbb{C}$, donc $\dim \mathcal{L}(0)=1$ et on a
$$\mbox{deg} K=g+\dim\mathcal{L}(K)-2.$$ Par ailleurs, en posant
$\mathcal{D}=0$ dans la formule (2), on obtient $$\dim
\mathcal{L}(0)-\dim\mathcal{L}(K)=\mbox{deg} 0-g+1,$$ d'où $ \dim
\mathcal{L}(K)=g$ et par conséquent $\mbox{deg} K=2g-2$. Ceci
achève la preuve du résultat annoncé. Pour terminer la preuve du
$2^{\mbox{ème}} cas$, on utilise ce résultat et la formule (7), on
obtient
$$
\dim \mathcal{L}(K-\mathcal{D})-\dim
\mathcal{L}(\mathcal{D})=-\mbox{deg}\mathcal{D}+g-1.
$$
$3^{\mbox{ème}} cas$ :
$\dim\mathcal{L}(\mathcal{D})=\dim\mathcal{L}(K-\mathcal{D})=0$.
Pour ce cas, on doit montrer que : $\mbox{deg} \mathcal{D}=g-1$.
Pour celà, considérons deux diviseurs positifs $\mathcal{D}_1$ et
$\mathcal{D}_2$ n'ayant aucun point en commun et posons
$\mathcal{D}\equiv \mathcal{D}_1-\mathcal{D}_2$. On a
$$\mbox{deg}\mathcal{D}=\mbox{deg} \mathcal{D}_1-\mbox{deg}
\mathcal{D}_2,$$ et
\begin{eqnarray}
\dim\mathcal{L}(\mathcal{D})&\geq& \mbox{deg}
\mathcal{D}_1-g+1,\nonumber\\
&=&\mbox{deg} \mathcal{D}+\mbox{deg} \mathcal{D}_2-g+1,\nonumber
\end{eqnarray}
i.e., $$\mbox{deg}\mathcal{D}_2-\dim
\mathcal{L}(\mathcal{D}_1)\leq \mbox{deg} \mathcal{D}+g-1.$$ Or
$$\mbox{deg}\mathcal{D}_2-\dim\mathcal{L}(\mathcal{D}_1)\geq 0,$$
car sinon il existe une fonction $f\in \mathcal{L}(\mathcal{D}_1)$
qui s'annule en tout point de $\mathcal{D}_2$, donc
$\mbox{deg}\mathcal{D}\leq g-1.$ En appliquant le même
raisonnement à $K-\mathcal{D}$, on obtient
$\mbox{deg}(K-\mathcal{D})\leq g-1.$ Comme $\mbox{deg}K=2g-2$
(voir (8)), alors $\mbox{deg}\mathcal{D}\geq g-1.$ Finalement,
$\mbox{deg}\mathcal{D}=g-1$, ce qui achève la démonstration du
théorème.

\begin{rem}
La formule (2) peut s'écrire sous la forme suivante :
$$\dim H^0(X,\mathcal{O}_\mathcal{D})-\dim H^1(X,\mathcal{O}_\mathcal{D})=\mbox{deg}
\mathcal{D}-g+1.$$ En introduisant la caractéristique
d'Euler-Poincaré :
$$
\chi(\mathcal{D})\equiv\dim \mathcal{L}(\mathcal{D})-\dim
\Omega(-\mathcal{D})=\dim H^0(X,\mathcal{O}_\mathcal{D})-\dim
H^1(X,\mathcal{O}_\mathcal{D}),
$$
pour un diviseur $\mathcal{D}$ sur une surface de Riemann $X$ de
genre $g$, le théorème de Riemann-Roch s'écrit
$$\chi(\mathcal{D})=\mbox{deg} \mathcal{D}-g+1.$$
\end{rem}

\begin{rem}
On sait que toute fonction holomorphe sur une surface de Riemann
compacte $X$ est constante. Une question se pose : Que se passe
t-il dans le cas des fonctions méromorphes? La réponse découle du
théorème de Riemann-Roch. Plus précisement, si $p$ est un point
quelconque de $X$, on peut trouver une fonction méromorphe non
constante, holomorphe sur $X\backslash\{p\}$ et ayant un pôle
d'ordre inférieur ou égal à $g+1$ en $p$. De même, on montre qu'il
existe sur $X$ des formes différentielles holomorphes non nulles,
qui s'annulent en au moins un point.
\end{rem}

\begin{rem}
Soient $(\omega_1,...,\omega_g)$ une base de $\Omega(X)$. Soit
$(U,\tau)$ une carte locale en $p\in X$ avec $\tau(p)=0$. Il
existe des fonctions $f_j$ holomorphes sur $U$ telles que :
$\omega_j=f_j(\tau)d\tau$. Le wronskien de $\omega_1,...,\omega_g$
est défini par le déterminant
$$W_\tau(\omega_1,...,\omega_g)\equiv W(f_1,...,f_g)=\det
(f^{(k-1)}_j)_{1\leq j,k\leq g}.$$ On dit que $p$ est un point de
Weierstrass si $W_\tau(\omega_1,...,\omega_g)$ s'annule. Dans le
cas où $p$ est un point de Weierstrass  alors on peut trouver une
fonction méromorphe sur $X$ ayant un pôle unique d'ordre inférieur
ou égal au genre $g$ au point $p$. Une autre application du
théorème de Riemann-Roch, permet de montrer l'existence d'une
suite de $g$ entiers : $1=n_1<n_2<...<2g,$ $g\geq1$, pour lesquels
il n'existe aucune fonction holomorphe sur $X\setminus {p}$, $p\in
X$, et ayant un pôle en $p$ d'ordre exactement $n_j$. On montre
que $p$ est un point de Weierstrass si et seulement si la suite
des $n_j$ est distincte de $\{1,2,...,g\}$.
\end{rem}

Soient $X$ et $Y$ deux surfaces de Riemann compactes connexes et
soit $f$ une application holomorphe non constante de $X$ dans $Y$.
Notons que $f$ est surjective, on dit que c'est un revêtement.
Pour tout point $p\in X$, il existe une carte $\varphi$ (resp.
$\psi$) de $X$ (resp. $Y$) centrée en $p$ (resp. $f(p)$) telles
que : $f_{\psi\circ\varphi}(\tau)=\tau^n$, où $n$ est un entier
strictement positif. L'entier $n-1$ s'appelle indice de
ramification de $f$ au point $p$ et on le note $V_p(f)$. Lorsque
$V_p(f)$ est strictement positif, alors on dit que $p$ est un
point de ramification (ou de branchement) de $f$. Une condition
nécessaire et suffisante pour que $p$ soit un point de
ramification de $f$ est que le rang de $f$ en $p$ soit nul.
L'image $J$ des points de ramifications de $f$ ainsi que que son
image réciproque $I$ sont fermés et discrets. La restriction de
$f$ à $X\setminus I$ est un revêtement de $Y\setminus J$ dont le
nombre de feuillets est le degré de l'application $f$ et on a
$$m\equiv \sum_{p\in f^{-1}(q)}(V_p(f)+1),\quad \forall q\in Y.$$

\begin{cor}
(Formule de Riemann-Hurwitz). Soient $X$ et $Y$ deux surfaces de
Riemann compactes de genre $g(X)$ et $g(Y)$ respectivement. Soit
$f$ une application holomorphe non constante de $X$ dans $Y$.
Alors
$$g(X)=m(g(Y)-1)+1+\frac{V}{2},$$
où $m$ est le degré de $f$ et $V$ est la somme des indices de
ramification de $f$ aux différents points de $X$.
\end{cor}
\textbf{Démonstration}: Soit $f:X\longrightarrow Y$, une
application holomorphe non constante de degré $m$. Désignons par
$\omega$ une forme différentielle méromorphe non nulle sur $Y$.
Soit $\tau$ (resp. $\upsilon$) un paramètre local sur $X$ (resp.
$Y$) et supposons que : $\upsilon=f(\tau)$. Soit
$\omega=g(\tau)d\tau$, une forme différentielle méromorphe sur $Y$
et soit $\eta=g(f(\tau))f'(\tau)d\tau$, une forme différentielle
sur $Y$. Nous allons voir que cette dernière est aussi méromorphe.
Notons que si on remplace $\tau$ par $\tau_1$, avec
$\tau=w(\tau_1)$, alors en terme de $\tau_1$ l'application $f$
s'écrit $\upsilon=(f\circ w)(\tau_1)$, et donc nous attribuons à
$\tau_1$ l'expression
$g(f(w(\tau_1)))f'(w(\tau_1))w'(\tau_1)d\tau_1$ ce qui montre que
$\eta$ est une forme différentielle méromophe. On peut supposer
que $\tau$ s'annule dans un voisinage de $p\in X$ et que
$\upsilon$ s'annule en $f(p)$. Dès lors,
$\upsilon=\tau^{V_f(p)+1}$ où $V_f(p)$ est l'indice de
ramification de $f$ au point $p$. Par conséquent,
$$\mbox{ord}_p\eta=(V_f(p)+1)\mbox{ord}_p\omega+V_f(p),$$ et
$$\sum_{p\in X}\mbox{ord}_p\eta=\sum_{p\in X}(V_f(p)+1)\mbox{ord}_{f(p)}\omega+V,$$
où $V=\sum_{p\in X}V_f(p)$. D'après la formule 8, on a
$$\sum_{p\in X}\mbox{ord}_p\eta=2g(X)-2,$$
et
\begin{eqnarray}
\sum_{p\in X}(V_f(p)+1)\mbox{ord}_{f(p)}\omega&=&\sum_{p\in X,
V_f(p)=0}\mbox{ord}_{f(p)}\omega,\nonumber\\
&=&\sum_{q\in Y}m . \mbox{ord}_q\omega,\nonumber\\
&=&m(2g(Y)-2).\nonumber
\end{eqnarray}
Par conséquent,
$$2g(X)-2=m(2g(Y)-2)+V,$$
ce qui achève la preuve du corollaire.

\section{Exemples}

Une des conséquences les plus intéressantes de la formule de
Riemann-Hurwitz est de donner un moyen efficace de calculer le
genre d'une surface de Riemann donnée.

\begin{Exmp}
Un cas particulier important est représenté par les courbes
hyperelliptiques $X$ de genre $g(X)$ d'équations
$$w^{2}=p_{n}(z)=\prod_{j=1}^n(z-z_j),$$
où $p_{n}(z)$ est un polynôme sans racines multiples, i.e., tous
les $z_j$ sont distincts. Notons que $$f:X\longrightarrow
Y=\mathbb{C}P^1=\mathbb{C}\cup \{\infty\},$$ est un revêtement
double ramifié le long des points $z_j$. Chaque $z_j$ est ramifié
d'indice $1$ et en outre le point à l'infini $\infty$ est ramifié
si et seulement si $n$ est impair. D'après la formule de
Riemann-Hurwitz, on a
\begin{eqnarray}
g(X)&=&m(g(Y)-1)+1+\frac{V}{2},\nonumber\\
&=&2(0-1)+1+\frac{1}{2}\sum_{p\in X}V_f(p),\nonumber\\
&=&E(\frac{n-1}{2}),\nonumber
\end{eqnarray}
où $E(\frac{n-1}{2})$ désigne la partie entière de
$(\frac{n-1}{2})$. Les courbes hyperelliptiques de genre $g$ sont
associées aux équations de la forme : $w^2=p_{2g+1}(z),$ ou
$w^2=p_{2g+2}(z),$ (selon que le point à l'infini $\infty$ est un
point de branchement ou non) avec $p_{2g+1}(z)$ et $p_{2g+2}(z)$
des polynômes sans racines multiples. Lorsque $g=1,$ on dit
courbes elliptiques.
\end{Exmp}

\begin{Exmp}
Déterminons le genre $g$ de la surface de Riemann $X$ associée à
l'équation :
$$F(w,z)=w^3+p_2(z)w^2+p_4(z)w+p_6(z)=0,$$
où $p_j(z)$ désigne un polynôme de degré $j$. On procède comme
suit : on a
\begin{eqnarray}
F(w,z)&=&w^3+az^2w^2+bz^4w+cz^6+\text{termes d'ordre
inférieur},\nonumber\\
&=&\prod_{j=1}^3(z+\alpha_jz^2)+\text{termes d'ordre
inférieur},\nonumber
\end{eqnarray}
Considérons $F$ comme un revêtement par rapport à $z $ et
cherchons ce qui ce passe quand $z \nearrow \infty .$ On a
\begin{eqnarray}
(w)_\infty&=&-2P-2Q-2R,\nonumber\\
(z)_\infty&=&-P-Q-R.\nonumber
\end{eqnarray}
Posons $t=\frac{1}{z}$, d'où
$$F(w,z)= \frac{1}{t^{6}}(t^6z^3+at^4z^2+bt^2z+c)+\cdots.$$
Ceci suggère le changement de cartes suivant: $$(w,z)\longmapsto
(\zeta=t^2w,t=\frac{1}{z }).$$ On a
\begin{eqnarray}
\frac{\partial F}{\partial w}&=&3w ^2+2p_2(z)w+p_4(z),\nonumber\\
&=&3w^2+2az^2w+bz^4+...,\nonumber\\
&=&\frac{3\zeta^2}{t^4}+\frac{2a\zeta}{t^4}+\frac{b}{t^4}+...\nonumber
\end{eqnarray}
La fonction $\frac{\partial F}{\partial w}$ étant méromorphe sur
la surface de Riemann $X$, alors Le nombre de zéros de cette
fonction coincide avec celui de ses pôles. Comme
\begin{eqnarray}
(\frac{\partial F}{\partial
w})_{P}&=&-4P,\nonumber\\
(\frac{\partial F}{\partial
w})_{Q}&=&-4Q,\nonumber\\
(\frac{\partial F}{\partial
w})_{R}&=&-4R,\nonumber\\
(\frac{\partial F}{\partial w})_{\infty}&=&-4(P+Q+R),\nonumber
\end{eqnarray}
alors le nombre de zéros de $\frac{\partial F}{\partial w }$ dans
la partie affine $X\setminus \{P,Q,R\}$ est égal à 8, et d'après
la formule de Riemann-Hurwitz, on a $g(X) =4.$
\end{Exmp}

\begin{Exmp}
Calculons le genre de la surface de Riemann $X$ associée au
polynôme :
$$w^4=z^4-1.$$
Ici, on a quatre feuillets. Les points de ramifications à distance
finie sont $1,-1,i$ et $-i$. On note que $z=\infty$ n'est pas un
point de ramification. L'indice de ramification étant égal à $12$,
alors d'après la formule de Riemann-Hurwitz, le genre de la
surface de Riemann en question est égal à $3$.
\end{Exmp}

\begin{Exmp}
Considérons la courbe de Fermat $X$ associée à l'équation :
$$w^n+z^n=1,\quad n\geq 2.$$
Ici on a un revêtement de degré $n$. Chaque racine
$n^{\mbox{ème}}$ de l'unité est ramifié d'indice $n-1$ tandis que
le point à l'infini $\infty$ n'est pas un point de ramification et
par conséquent $$g(X)=\frac{(n-1)(n-2)}{2}.$$ L'équation de Fermat
: $$U^n+V^n=W^n,$$  (avec $w=\frac{U}{Z}$, $z=\frac{V}{Z}$) étant
de genre $\geq 1$ pour $n\geq 3$, elle n'admet donc qu'un nombre
fini de solutions. Ce fut une des pistes utilisées récemment par
A. Wiles pour prouver le grand théorème de Fermat : pour $n\geq 3$
cette équation n'a pas de solution non triviale.
\end{Exmp}

\begin{Exmp}
Déterminons le genre de la surface de Riemann $X$ associée au
polynôme :
$$(w^2-1)((w^2-1)z^4-p(z))+c=0,$$
où $$p(z)=az^2-2bz-1,$$ et $a,b,c$ sont des constantes non nulles.
Cette surface a été obtenue pour la première fois par l'auteur [6]
en 1988, lors de l'étude du célèbre problème de Kowalewski
concernant la rotation d'un corps solide autour d'un point fixe.
Cette surface a permit de bien comprendre la géométrie liée à ce
problème et a été utilisée par la suite pour élucider d'autres
questions (voir par exemple [1]) concernant le flot géodésique sur
le groupe des rotations $SO(4)$ ainsi que celui du système de
Hénon-Heiles. Le calcul du genre de $X$ n'est pas immédiat. Notons
tout d'abord que l'application $$\sigma:X\longrightarrow X,
(w,z)\longmapsto (-w,z),$$ est une involution (automorphisme
d'ordre deux) sur $X$. Le quotient $Y=X/\sigma$ de $X$ par
l'involution $\sigma$ est une courbe elliptique définie par
$$u^2=p^2(z)-4cz^4.$$
La surface de Riemann $X$ est un revêtement double ramifié le long
de la courbe $Y$ : $$X\longrightarrow Y, \quad(w,u,z)\longmapsto
(u,z),$$
$$
X : \left\{\begin{array}{rl} w^2&=\frac{2z^4+p(z)+u}{2z^4}\\
u^2&=p^2(z)-4cz^4
\end{array}\right.
$$
Pour $z$ suffisament petit, on a
$$w^2=\frac{2z^4+p(z)+\sqrt{p^2(z)-4cz^4}}{2z^4}=1-c+\circ (z),$$
et
$$w^2=\frac{2z^4+p(z)-\sqrt{p^2(z)-4cz^4}}{2z^4}=\frac{1}{z^4}(-1+\circ (z)).$$
Au voisinage de $z=\infty$, on a
$$2(w^2-1)z^2=a\pm \sqrt{a^2-4c}+\circ (z).$$
La surface $X$ possède quatre points à l'infini $p_j$ $(1\leq
j\leq 4)$ et quatre points de ramifications $$q_j\equiv (w=0,
u=-2z^4-p(z), z^4+p(z)+c=0),\quad 1\leq j\leq 4,$$ sur la courbe
elliptique $Y$. Dès lors, la structure des diviseurs de $w$ et $z$
sur $X$ est
\begin{eqnarray}
(w)&=&\sum_{1\leq j\leq 4}q_j-\sum_{1\leq j\leq 4}p_j,\nonumber\\
(w)&=&\mbox{quatre zéros}-\sum_{1\leq j\leq 4}p_j.\nonumber
\end{eqnarray}
Finalement, en appliquant la formule de Riemann-Hurwitz, on
obtient
$$g(X)=m(g(Y)-1)+1+\frac{V}{2}=2(1-1)+1+\frac{4}{2}=3.$$
\end{Exmp}

\end{document}